\documentclass[12pt]{article}

\usepackage{amsmath}
\usepackage{amsfonts}
\usepackage{amsthm}

\theoremstyle{remark}

\numberwithin{equation}{section}

\begin{document}

\begin{center}
{\Large {\bf A Property of K\"{a}hler-Ricci Solitons \ on Complete Complex Surfaces}} \vskip 5mm Bing-Long Chen and Xi-Ping Zhu

\smallskip 
Department of Mathematics,\\Zhongshan University,\\
Guangzhou 510275, P. R. China\\ 
\medskip
\end{center}

\baselineskip=20pt
\section*{1. Introduction}

\setcounter{section}{1} \setcounter{equation}{0} \qquad This paper is
concerned with Ricci flow on complete K\"{a}hler manifolds. Let $M$ be a complete complex manifolds with K\"{a}hler metric $g_{\alpha \overline{\beta}}$. The Ricci flow is the following evolution equation on the metric
\begin{equation}
\label{1.1}\left\{ 
\begin{array}{ll}
\bigbreak \displaystyle \frac{\partial g_{\alpha \overline{\beta }}(x,t)}{%
\partial t}=-R_{\alpha \overline{\beta }}(x,t)\  & \qquad x\in M\ \quad
t>0\ , \\ 
g_{\alpha \overline{\beta }}(x,0)=g_{\alpha \overline{\beta }}(x)\  & 
\qquad x\in M\ ,
\end{array}
\right. 
\end{equation}
where $R_{\alpha \overline{\beta }}(x,t)$ denotes the Ricci curvature tensor
of the metric $g_{\alpha \overline{\beta }}(x,t)$.

One of the main problems in differential geometry is to find canonical structure on manifolds. The Ricci flow introduced by Hamilton [8] is an useful tool to approach such problems. For examples, Hamilton [10] and Chow [7] used the convergence of the Ricci flow to characterize the complex structures on compact Riemann surfaces, Hamilton [8] used the Ricci flow to classify compact three-manifolds with positive Ricci curvature, and the authors [4] recently used the Ricci flow to get steinness for a class of complete noncompact K\"{a}hler manifolds. By a direct computation one can see that the scalar curvature $R(x,t)$ of $g_{\alpha \overline{\beta}}(x,t)$ satisfies the equation 
$$\frac \partial {\partial t}R = \bigtriangleup R+2\left|R_{\alpha \overline{\beta}}\right| ^2.$$
This is a nonlinear heat equation with superlinear growth. It is clear that the scalar curvature must generally blow up in finite time. In the other words the Ricci flow will develop singularities in finite time. Thus it is important to consider what kind of singularities might form.

In [12], Hamilton divided the solutions of the Ricci flow into three types and showed that one can dilate the solutions around the singularities to get in limit the models of the corresponding type. Then Cao [2] showed that if in addition the limits of Type II or Type III have nonnegative holomorphic bisectional curvature, the singularity models are the solutions to (1.1) which move by one-parameter group of biholomorphisms or also expand by a factor at the same time.

A solution to the Ricci flow (1.1) which moves by a one-parameter group of biholomorphisms is called a K\"{a}hler-Ricci soliton. The equation for a K\"{a}hler metric to move by a biholomorphism in the direction of a holomorphic vector field $V$ is that the Ricci term $R_{\alpha \overline{\beta}}$ is the Lie deivative of the metric $g_{\alpha \overline{\beta}}$ in the direction of the vector field $V$; thus
\begin{equation}
\label{1.2}\left\{ 
\begin{array}{ll}
\bigbreak \displaystyle  R_{\alpha \overline{\beta }}=g_{\alpha \overline{\gamma }}\bigtriangledown_{\overline{\beta}}V^{\overline{\gamma}} + g_{\gamma \overline{\beta }}\bigtriangledown_{\alpha}V^{\gamma}, \\ 
\bigtriangledown_{\alpha}V^{\overline{\beta}} = \bigtriangledown_{\overline{\beta}}V^{\alpha} = 0.
\end{array}
\right. 
\end{equation}
From [12], any K\"{a}hler-Ricci soliton must be noncompact. Recall that a complete noncompact complex $n$-dimensional K\"{a}hler manifold $(M,g_{\alpha \overline{\beta }})$ is of maximal volume growth if there exist some fixed point $x_{0} \in M$ and some positive constant $C_{0}$ such that 
$$
\mbox{Vol}\,(B(x_0,r))\geq C_0r^{2n}\  
 \qquad \qquad \mbox{for all}\qquad 0\leq r<+\infty .
$$
where $\mbox{Vol}\,(B(x_0,r))$ denotes the volume of the geodesic ball $B(x_0,r)$ centered at $x_0$ with radius $r$.
The main result of this paper is the following property for K\"{a}hler-Ricci solitons.\vskip  3mm

{\bf \underline{Theorem}} \ \ \ Let $(M,g_{\alpha \overline{\beta }})$ be a nonflat K\"{a}hler-Ricci soliton on a complex two-dimensional K\"{a}hler manifold with bounded and nonnegative holomorphic bisectional curvature. Then $(M,g_{\alpha \overline{\beta }})$ can not be of maximal volume growth. \vskip  3mm

We conjecture this theorem holds for all dimensions. In the forthcoming paper [5], we will use this result to obtain a uniformization theorem for complete noncompact K\"{a}hler surfaces with positive holomorphic bisectinal curvature.

This paper contains three sections. In the next section we will show that the curvature operator of the K\"{a}hler-Ricci soliton is nonnegative definite everywhere. In the last section we will use a blow down argument to prove the theorem.

We are grateful to H.D. Cao for many helpful discusions. This work was partially supported by the Foundation for Outerstanding Young Scholars of China.

\section*{2. Nonnegativity of Curvature Operator}

\setcounter{section}{2} \setcounter{equation}{0}\qquad  The solution to the Ricci flow (1.1) corresponding to the K\"{a}hler-Ricci soliton exists for $-\infty < t <+\infty$ and is obtained by flowing along the holomorphic vector field $V$ in (1.2). The following lemma shows that the nonnegativity of holomorphic bisectional curvature of the K\"{a}hler-Ricci soliton on a complex surface implies the nonnegativity of its curvature operator.\vskip  3mm

{\bf \underline{Lemma}} \ \ \ Let $g_{\alpha \overline{\beta }}(x,t)$ be a complete solution to the Ricci flow (1.1) on a complex two-dimensional K\"{a}hler manifold $M$ for all $t \in (-\infty, +\infty)$. Suppose its holomorphic bisectional curvature is bounded and nonnegative everywhere. Then the curvature operator of the metric $g_{\alpha \overline{\beta }}(x,t)$ is nonnegative definite everywhere on $M\times (-\infty, +\infty)$. \vskip  3mm

{\bf \underline{Proof.}} \quad 
Choose a local orthonormal coframe 
$\{\omega _1,\omega _2,\omega _3,\omega _4\}$ 
on an open set 
$U\subset M$ so that 
$\omega _1+\sqrt{-1}\omega _2$ and $\omega_3+\sqrt{-1}\omega _4$ 
are $(1,0)$-forms over $U$. Then the self--dual forms 
$$
\varphi _1=\omega _1\land \omega _2+\omega _3\land \omega _4, 
\,\,\, \varphi_2=\omega _2\land \omega _3+\omega _1\land \omega _4, 
\,\,\, \varphi _3=\omega_3\land \omega _1+\omega _2\land \omega _4
$$ 
and the anti--self--dual forms 
$$
\psi_1=\omega _1\land \omega _2-\omega _3\land \omega _4, 
\,\,\, \psi _2=\omega_2\land \omega _3-\omega _1\land \omega _4, 
\,\,\, \psi _3=\omega _3\land \omega_1-\omega _2\land \omega _4
$$ 
form a basis of the space of $2$-forms over $U$.
In particular, $\varphi _1,\psi _1,\psi _2$ and $\psi _3$ 
give a basis for the space of $(1,1)$-forms over $U$.

On a complex two-dimensional K\"ahler manifold, it is well known that its curvature 
operator has image in the holonomy algebra $u(2) \,(\subset so(4))$ 
spanned by $(1,1)$-forms. Thus, the curvature operator $R_m$ 
in the basis $\{\varphi_1,\varphi _2,\varphi _3,\psi _1,\psi _2,\psi _3\}$ 
has the following form,
$$
R_m=\left( 
\begin{array}{cc}
\begin{array}{ccc}
a & 0 & 0 \\ 
0 & 0 & 0 \\ 
0 & 0 & 0 
\end{array}
& 
\begin{array}{ccc}
b_1 & b_2 & b_3 \\ 
0 & 0 & 0 \\ 
0 & 0 & 0 
\end{array}
\\ 
\begin{array}{ccc}
b_1 & 0 & 0 \\ 
b_2 & 0 & 0 \\ 
b_3 & 0 & 0 
\end{array}
& A 
\end{array}
\right) \ 
$$
where $A$ is a $3\times 3$ symmetric matrix.

Let $V$ be a real tangent vector of the complex surface 
$M$. Denote by $J$ the complex structure. It is clear that the complex 
$2$--plane $V\wedge JV$ is dual to $(1,1)$-form 
$u\varphi _1+v_1\psi _1+v_2\psi _2+v_3\psi _3$ satisfying the
decomposability condition $u^2=v_1^2+v_2^2+v_3^2$. Then after
normalizing $u$ to 1 by scaling, we see that the holomorphic 
bisectional curvature is nonnegative if and only if 
\begin{equation}
\label{2.1}
a+b\cdot v+b\cdot w+^{\ t}\negthinspace {}vAw\geq 0\ , 
\end{equation}
for any unit vectors $v=(v_1,v_2,v_3)$ and $w=(w_1,w_2,w_3)$ 
in ${\bf R^3}$, where $b$ is the vector $(b_1,b_2,b_3)$ 
in $R_m$.

Denote by $a_1\leq a_2\leq a_3$ the eigenvalues of $A$. 
Recall that $\mbox{tr}\,A=a $ by the Bianchi identity, 
so if we choose $v$ to be the eigenvector of $A $ with
eigenvalue $a_3$ and choose $w=-v$, (2.1) gives 
\begin{equation}
\label{2.2}a_1+a_2\geq 0\ .
\end{equation}
In particular, we have $a_2\geq 0$.

To proceed further, we need to adapt the Hamilton's maximum principle for tensors.
Let
$$
\left( a_i\right) _{\min }(t)=\inf\limits_{x\in M}a_i(x,t)\
,\qquad i=1,2,3\ 
$$
and
$$
K=\sup \limits_{(x,t)\in M\times (-\infty ,+\infty)}\left|
Rm(x,t)\right|.
$$
By assumption, the solution
$g_{\alpha \overline{\beta }}(x,t), (-\infty < t < +\infty),$
has bounded holomorphic bisectional curvature, hence $K$ is finite.
Thus, by the derivative estimate of Shi [13] 
(see also Theorem 7.1 in [12]), the all 
derivatives of the curvature are also uniformly bounded. 
In particular, we can use the maximum
principle of Cheng--Yau (see Proposition 1.6 in [6]) and then, as
observed in [12], this implies that the maximum principle of 
Hamilton in [9] actually works for the evolution equations 
of the curvature of $\widetilde{g}_{\alpha \overline{\beta }}(x,t)$ 
on the complete noncompact manifold $M.$ Thus, from [9],
we obtain
\begin{eqnarray}
\frac{d\left( a_1\right) _{\min }}{dt}&\geq&\left( \left( a_1\right)_
{\min }\right) ^2+2\left( a_2\right) _{\min }\cdot 
\left( a_3\right) _{\min } \nonumber \\
& \geq & 3\left( \left( a_1\right) _{\min }\right) ^2 \nonumber
\end{eqnarray}
by (2.2). Then, for fixed $t_0\in (-\infty ,0)$ and 
$t>t_0,$
\begin{eqnarray}
\left( a_1\right) _{\min }(t)&\geq&\frac 1{\left( a_1\right)_
{\min }^{-1}\left( t_0\right) -3\left( t-t_0\right) } \nonumber \\
& \geq & \frac 1{-K^{-1}-3\left( t-t_0\right) }\ . \nonumber
\end{eqnarray}
Letting $t_0\rightarrow -\infty $, we get 
\begin{equation}
\label{2.3}
a_1\geq 0\ ,\qquad \mbox{for all} \;\;\; (x,t) \in 
M \times (-\infty ,0]\,
\end{equation}
i.e. $A$ is nonnegative definity.

Finally, to prove the nonnegativity of the curvature operator 
$R_m$, we recall its corresponding ODE from [9],
$$
\frac {dR_m}{dt}={R_m}^2+\left( 
\begin{array}{cc}
0 & 0 \\ 
0 & A^{\#} 
\end{array}
\right) \ , 
$$
where $A^{\#}$ is the adjoint matrix of $A.$

Let $m_1$ be the smallest eigenvalue of the curvature operator $R_m$. 
Exactly as above, by using the maximum principle of Hamilton, we have
$$
\frac{d\left( m_1\right) _{\min }}{dt}\geq \left( m_1\right) _{\min }^2 
$$
where $\left( m_1\right) _{\min }(t)=\inf \limits_{x\in M}m_1(x,t)$. Therefore, by the same reasoning in the derivation of 
(2.3), we have
\begin{equation}
\label{3.7}
m_1\geq 0\ ,\qquad \mbox{for all} \;\; (x,t)\in M\times
(-\infty , +\infty). 
\end{equation}
So $R_m$ is nonnegative definite everywhere and the proof of the lemma is completed.
\hfill
$Q. E. C.$

\section*{3. Blow Down Argument}

\setcounter{section}{3} \setcounter{equation}{0}\ ~\quad  This section is devoted to the proof of the main result.

\vskip 3mm

{{\underline{{\bf The proof of the theorem.}}}}

Let $(M,g_{\alpha \overline{\beta }}(x))$ be a nonflat K\"{a}hler-Ricci soliton on a complex two-dimensional K\"{a}hler manifold with bounded and nonnegative holomorphic bisectional curvature. Denote $g_{\alpha \overline{\beta }}(x,t), t \in (-\infty, +\infty),$ to be the solution to the Ricci flow (1.1) corresponding to the K\"{a}hler-Ricci soliton $g_{\alpha \overline{\beta }}(x)$.
The metric $g_{\alpha \overline{\beta }}(x,t)$ is obtained from the pullback of $g_{\alpha \overline{\beta }}(x)$ by a one-parameter group of biholomorphism. Of course, for each $t \in (-\infty, +\infty),g_{\alpha \overline{\beta }}(x,t)$ also has bounded and nonnegative holomorphic bisectional curvature.

We prove by contradiction. Thus, suppose $g_{\alpha \overline{\beta }}(x)$ is of maximal volume growth. This is, there exist some point $x_0 \in M$ and some positive constant $C_0$ such that  
\begin{equation}
\mbox{Vol}\left(B(x_0,r)\right) \geq C_0r^4\ 
\qquad \mbox{for all} \quad 0\leq r< +\infty . 
\end{equation}
Note that $g_{\alpha \overline{\beta }}(x)$ has nonnegative Ricci curvature. Then it follows from the standard volume comparison that the above inequality actually holds for all $x_0 \in M$. Hence
\begin{equation}
\label{3.1}
\mbox{Vol}_t\left(B_t(x,r)\right) \geq C_0r^4\ 
\qquad \mbox{for all} \quad 0\leq r< +\infty \quad \mbox{and} 
\quad x\in M\ , 
\end{equation}
where $\mbox{Vol}_t\left(B_t(x,r)\right)$ denotes 
the volume of the geodesic ball $B_t(x,r)$ of radius $r$ with 
center at $x$ with respect to the metric 
$g_{\alpha \overline{\beta }}(x,t).$

We first observe that
\begin{equation}
\label{3.3}
\limsup \limits_{d(x,x_0)\rightarrow +\infty } 
R(x)d^2(x,x_0)=+\infty \  
\end{equation}
where $d(x,x_0)$ the geodesic distance between two points 
$x,x_0\in M$ with respect to the metric $g_{\alpha \overline{\beta }}(x)$.
 
In fact, suppose not, thus the curvature of the initial metric 
$g_{\alpha \overline{\beta }}(x,0)=g_{\alpha \overline{\beta }}(x)$ in (1.1)
has quadratic decay. By applying a result of Shi 
(see Theorem 8.2 in [14]), we see that
\begin{equation}
\lim \limits_{t\rightarrow +\infty }\sup \left\{ \left. 
R(x,t)\right| \ x\in M\right\} =0\ . 
\end{equation}
where $R(x,t)$ is the scalar curvature of the solution $g_{\alpha \overline{\beta }}(x,t)$. But $g_{\alpha \overline{\beta }}(x,t)$ is a nonflat K\"{a}hler-Ricci soliton. This is impossible.

With the estimate (3.3), we can then apply a lemma of Hamilton 
(Lemma 22.2 in [12]) to find a sequence of points $x_{j}$, a sequence of radii $r_j$ and a sequence of positive numbers 
$\delta _j$, $j=1,2,\cdots$, with $\delta_j\rightarrow 0$ such that

\begin{enumerate}
\item[{(a)}]  $R(x,0)\leq (1+\delta _j)R(x_j,0)$ 
for all $x$ in the ball $B(x_j,r_j)$ of radius $r_j$ 
centered at $x_j$ with respect to the metric 
$g_{\alpha \overline{\beta }}$;

\item[{(b)}]  $r_j^2 R(x_j,0)\rightarrow +\infty $;

\item[{(c)}]  if $s_j=d(x_j,x_0)$, then $\lambda
_j=s_j/r_j\rightarrow +\infty $;

\item[{(d)}]  the balls $B(x_j,r_j)$ are disjoint.
\end{enumerate}

We have shown in the previous section that the metric $g_{\alpha \overline{\beta }}$ has nonnegative definite curvature operator. In particular, the sectional curvature is nonnegative. Denote the minimum of the sectional curvature of the metric $g_{\alpha \overline{\beta }}$ at $x_j$ by $\nu_j$. We claim that the following holds 
\begin{equation}
\label{3.5}
\varepsilon _j=\frac{\nu_j}{R(x_j,0)}\rightarrow 0\
\qquad \mbox{as} \quad j\rightarrow +\infty \ . 
\end{equation}

Suppose not, there exists a subsequence $j_k\rightarrow +\infty $ 
and some positive number $\varepsilon >0$ such that 
\begin{equation}
\label{3.6}
\varepsilon _{j_k}=\frac{\nu_{j_k}}{R(x_{j_k},0)}\geq
\varepsilon \ \qquad \mbox{for all} \quad k=1,2,\cdots \ . 
\end{equation}

Since the solution $g_{\alpha \overline{\beta }}(x,t)$
exits for $-\infty < t < +\infty$, it follows from the Li-Yau type differential inequality of Cao [1] that
the scalar curvature $R(x,t)$ is pointwisely nondecreasing in
time. Then, by using the local derivative estimate of Shi [13] 
(or see Theorem 13.1 in [12]) and (a), (b), we have
\begin{eqnarray}
\label{3.7}
\sup \limits_{x\in B(x_{j_k},r_{j_k})}\left| 
\nabla R_m(x,0)\right| ^2&\leq&C_1R^2(x_j,0)
\left( \frac 1{r_{j_k}^2}+R(x_j,0)\right) \nonumber \\
&\leq&2C_1R^3(x_j,0)\ ,
\end{eqnarray}
where $C_1$ is a 
positive constant depending only on the dimension.

For any $x\in B(x_{j_k},r_{j_k})$, we obtain 
from (3.6) and (3.7) that the minimum of the sectional curvature $\nu(x)$ of $g_{\alpha \overline{\beta}}$ at $x$, satisfies
\begin{eqnarray}
\label{3.8}
\nu(x)&\geq&\nu_{j_k}-\sqrt{2C_1}R^{3/2}(x_{j_k},0)d(x,x_{j_k})\nonumber\\
&\geq&R(x_{j_k},0)\left( \varepsilon -\sqrt{2C_1}\cdot \sqrt{R(x_{j_k},0)}\cdot d(x,x_{j_k})\right)\nonumber\\
&\geq&\frac {\varepsilon }2 R(x_{j_k},0)\ 
\end{eqnarray}
if
$$
d(x,x_{j_k})\leq \frac{\varepsilon}{2\sqrt{2C_1}\cdot \sqrt{ 
R(x_{j_k},0)}}\ . 
$$
Thus, from (a) and (3.8), there exists $k_0>0$ such that
for any $k\geq k_0$ and 
$$
x \in B(x_{j_k},\frac{\varepsilon}{2\sqrt{2C_1}
\cdot \sqrt{R(x_{j_k},0)}}),
$$
we have
\begin{equation}
\frac {\varepsilon} 2 R(x_{j_k},0) \leq 
\mbox{the sectional curvature at $x$} \leq 
2R(x_{j_k},0). 
\end{equation}

Therefore the balls $B(x_{j_k},\frac{\varepsilon}{2\sqrt{2C_1}
\cdot \sqrt{R(x_{j_k},0)}}), k_0 \leq k < +\infty,$ are a family of disjont remote curvature $\beta$-bumps for some $\beta > 0$ in the sense of Hamilton [12]. But this contradicts with the finite bumps theorem of Hamilton [12]. So we have proven the claim (3.5).

Next we blow down the K\"{a}hler-Ricci soliton $g_{\alpha \overline{\beta }}(x,t))$ along the points $x_j$. For the above chosen $x_j$, 
$r_j$ and $\delta _j$, let $x_j$ be the new origin, dilate the space by
a factor $\lambda _j$ so that $R(x_j,0)$ become $1$ at the
origin at $t=0$, and dilate in time by $\lambda _j^2$ so that it is still a
solution to the Ricci flow. The balls $B(x_j,r_j)$ are dilated
to the balls centered at the origin of radii $\widetilde{r}_j=r_j^2 
R(x_j,0)\rightarrow +\infty $ ( by (b) ). Since the scalar
curvature of $g_{\alpha \overline{\beta }}(x,t)$ is pointwise
nondecreasing in time by the Cao's inequality [1], the curvature bounds on $B(x_j,r_j)$ also give bounds for previous times in these balls. 

On the other hand, as shown in [4], the combination of the maximal volume growth estimate (3.1) and the local
injectivity radius estimate of Cheeger, Gromov and Taylor [3] implies
that
$$
\mbox{inj}_{M}\left( x_j,g_{\alpha 
\overline{\beta }}\right) 
\geq \frac \alpha {\sqrt{R(x_j,0)}}\ ,
$$
for some positive constant $\alpha $ independent of $j.$

So we have everything to take a limit for the dilated solutions. By applying
the compactness theorem in [11] and combining (3.2), (3.5), (a) and (b), we obtain a complete noncompact solution
$( \widetilde{M},\widetilde{g}_{\alpha \overline{\beta }}(x,t))$ of the Ricci flow for $t\in
(-\infty ,0]$ such that

\begin{enumerate}
\item[{(e)}]  the curvature operator is still nonnegative;

\item[{(f)}]  the scalar curvature $\widetilde{R}(x,t)\leq 1$, for all $x\in \widetilde{M}$, $
t\in (-\infty ,0]$, and $\widetilde{R}(0,0)=1$;

\item[{(g)}]  the volume of the geodesic ball to metric $\widetilde{g}_{\alpha \overline{\beta }}(x,t)$ still satisfies
 $$\mbox{Vol}_t\left( \widetilde{B}_t(x,r)\right) \geq C_0r^4$$ 
for all $x\in \widetilde{M}$, $0\leq r < +\infty $;

\item[{(h)}]  there exists a complex $2$--plane at the origin so that at $t=0$, the corresponding sectional curvature vanishes.
\end{enumerate}

If we consider the universal covering of $\widetilde{M}$, the induced metric
of $\widetilde{g}_{\alpha \overline{\beta }}(x,t)$ on the universal
covering is clearly still a solution to the Ricci flow and satisfies all 
of above (e), (f), (g), (h). Thus, without loss of generality, we may assume 
that $\widetilde{M}$ is simply connected.

By using the strong maximum principle on the evolution equation of the
curvature operator of $\widetilde{g}_{\alpha \overline{\beta }}(x,t)$
as in [9] (see Theorem 8.3 of [9]), we know that there
exists a constant $K>0$ such that on the time interval $-\infty <t<-K$, the
image of the curvature operator of $(\widetilde{M},\widetilde{g}_{\alpha 
\overline{\beta }}(x,t))$ is a fixed Lie subalgebra of $so(4)$ of
constant rank on $\widetilde{M}$. Because $\widetilde{M}$ is K\"ahler, the
possibilities are limited to $u(2)$, $so(2)\times so(2)$ or $so(2).$

In the case $u(2)$, the sectional curvature is strictly positive. Thus, this
case is ruled out by (h). In the cases $so(2)\times so(2)$ or $so(2)$,
according to [9], the simply connected manifold $\widetilde{M}$
splits as a product $\widetilde{M}=\Sigma _1\times \Sigma _2$, where $\Sigma
_1$ and $\Sigma _2$ are two Riemann surfaces with nonnegative curvature (by
(e)), and at least one of them, say $\Sigma _1$, has positive curvature (by (f)).

Denote by $\widetilde{g}_{\alpha \overline{\beta }}^{(1)}(x,t)$ the
corresponding metric on $\Sigma _1$. Clearly, it follows from (g) and
standard volume comparison that for any $x\in \Sigma _1$, $t\in (-\infty
,-K) $, we have 
\begin{equation}
\label{3.10}
\mbox{Vol}B_{\Sigma_1}(x,r) \geq C_2 r^2 \qquad \mbox{for} \;\; 
0 \leq r  < +\infty
\end{equation}
where both the geodesic ball $B_{\Sigma_1}(x,r)$ 
and the volume are taken with respect to the metric 
$\widetilde{g}_{\alpha \overline{\beta }}^{(1)}(x,t)$ on $\Sigma_1$,
$C_2$ is a positive constant depending only on $C_0.$
Also as the curvature of 
$\widetilde{g}_{\alpha \overline{\beta }}^{(1)}(x,t)$ 
is positive, it follows from Cohn--Vossen inequality that 
\begin{equation}
\label{3.11}
\int_{\Sigma _1}\widetilde{R}^{(1)}(x,t)d\sigma _t\leq 8\pi \ , 
\end{equation}
where $\widetilde{R}^{(1)}(x,t)$ is the scalar curvature of $(\Sigma _1, 
\widetilde{g}_{\alpha \overline{\beta }}^{(1)}(x,t))$ and $d\sigma _t$ is
the volume element of the metric $\widetilde{g}_{\alpha \overline{\beta }%
}^{(1)}(x,t).$

Now, the metric 
$\widetilde{g}_{\alpha \overline{\beta }}^{(1)}(x,t)$ 
is a solution to the Ricci flow on the Riemann surface $\Sigma _1$ 
over the time interval $(-\infty ,-K)$. 
The estimates (3.10) and (3.11) imply that for each $t\in (-\infty ,-K)$,
the curvature of $\widetilde{g}_{\alpha \overline{\beta }}^{(1)}(x,t)$ has
quadratic decay in the average sense of Shi [14] and then the a priori
estimate of Shi (see Theorem 8.2 in [14]) implies that the solution 
$\widetilde{g}_{\alpha \overline{\beta }}^{(1)}(x,t)$ exists for all 
$t\in (-\infty ,+\infty )$ and satisfies
\begin{equation}
\label{3.12}
\lim \limits_{t\rightarrow +\infty }\sup \left\{ \left. 
\widetilde{R}^{(1)}(x,t)\right| \ x\in \Sigma _1\right\} =0.
\end{equation}
Again, by the inequality of Cao [1], 
we conclude that
$$
\widetilde{R}^{(1)}(x,t)\equiv 0\qquad \mbox{on}\quad 
\Sigma _1\times (-\infty,+\infty )\ . 
$$
This contradicts with the fact that $(\Sigma _1,\widetilde{g}_{\alpha 
\overline{\beta }}(x,t))$ has positive curvature for $t<-K.$ Hence we have completed the 
proof of the theorem.
\hfill $Q. E. C.$

\end{document}